\documentclass{amsart}
\usepackage{amsthm}
\usepackage{amsmath}
\usepackage{amssymb}
\usepackage{booktabs}

\newtheorem{theorem}{Theorem}[section] % 1st argument is your name for it
\title{Zeroes of partial sums of the zeta-function}
\author{D.J. Platt}
\address{Heilbronn Institute for Mathematical Research, University of Bristol, Bristol, UK}
\email{dave.platt@bris.ac.uk}

\author{T.S. Trudgian}
\address{Mathematical Sciences Institute, The Australian National University,
 ACT 0200, Australia}
\email{timothy.trudgian@anu.edu.au}

\thanks{We are grateful to Andrew Booker with whom we discussed several aspects of this problem and to the referees for their careful reading of our manuscript. The second author was supported by Australian Research Council DECRA Grant DE120100173.}

\begin{document}
\maketitle

\begin{abstract}
%\noindent
This article considers the positive integers $N$ for which $\zeta_{N}(s) = \sum_{n=1}^{N} n^{-s}$ has zeroes in the half-plane $\Re(s)>1$. Building on earlier results, we show that there are no zeroes for $1\leq N\leq 18$ and for $N=20, 21, 28$. For all other $N$ there are infinitely many such zeroes.
\end{abstract}

\section{Introduction}
The Riemann zeta-function is defined as $\zeta(s) = \sum_{n=1}^{\infty} n^{-s}$ for $\Re(s) >1$. Throughout this article we write the complex variable $s$ as $s = \sigma + it$ with $\sigma$ and $t$ real, and consider $N$ to be a natural number. Truncation of the zeta-function gives the partial sum $\zeta_{N}(s) = 1 + 2^{-s} + \cdots + N^{-s}$. One may study these partial sums in the hope of deducing some information about $\zeta(s)$. For a comprehensive treatment of these ideas we refer the reader to \cite{GonekZero} and \cite{GonekMontyZero}. 

Tur\'{a}n \cite{Turan} showed that the Riemann hypothesis would follow if for all $N$ sufficiently large $\zeta_{N}(s)$ had no zero in $\sigma>1$. Let $\psi_{N}$ be the supremum over all values of $\sigma$ for which $\zeta_{N}(s)= 0$. Montgomery \cite{MontyU} showed that for all $N$ sufficiently large,
\begin{equation*}
\psi_{N} = 1 + \left( \frac{4}{\pi} - 1 - o(1)\right) \frac{\log\log N}{\log N},
\end{equation*}
where the constant $4/\pi - 1$ is best possible.
Therefore for $N$ sufficiently large, $\zeta_{N}(s)$ has zeroes in $\sigma>1$. 
 
Monach \cite{Monach} made this explicit: for all $N>30$ there are zeroes in $\sigma>1$. His proof was in two parts: an analytic argument for $N\geq 549,798$ and a computational proof for $30<N<549,798$. The latter proof is contained in \cite[Lem.\ 3.14, pp.\ 134-135]{Monach}. Monach's work can be combined with the results of Tur\'{a}n and Spira to give the following table.

\begin{table}[ht]
\caption{Zeroes of $\zeta_{N}(s)$ in $\sigma>1$ for various values of $N$.}
\label{table1}
\centering
\begin{tabular}{c c }
\hline\hline
Range of $N$ & Are there zeroes of $\zeta_{N}(s)$ in $\Re(s)>1$? \\[0.5ex]\hline
$1-5$ & No, \cite[pp.\ 7-8]{Turan}\\
$6-9$ & No, \cite[p.\ 550]{SpiraSections1} and \cite[Table II, \S 4]{SpiraSections2}\\
$19$ & Yes, \cite[Table III, \S 4]{SpiraSections2}\\
$22-27$ & Yes, \cite[Table III, \S 4]{SpiraSections2}\\
$29-50$ & Yes, \cite[Table III, \S 4]{SpiraSections2}\\
$\geq 51$ & Yes, \cite[Thm.\ 3.8]{Monach}\\ 
    \hline\hline
  \end{tabular}
\end{table}
Indeed, van de Lune and te Riele \cite{vantR} actually computed some zeroes of $\zeta_{N}(s)$ for $N=19, 22-27, 29-35, 37-41, 47.$ Adapting Bohr's theorems on values assumed by Dirichlet series, Spira \cite[Thm.\ 3]{Spiraset} (see also, \cite[p.\ 163]{SpiraSections2}) showed that if $\zeta_{N}(s)$ has one zero in $\sigma>1$ then it has infinitely many such zeroes. 

Therefore, all that remains is to investigate whether, for,
\begin{equation}\label{big}
N\in\{10,11,12,13, 14, 15, 16, 17, 18, 20, 21, 28\},
\end{equation}
$\zeta_{N}(s)$ has zeroes in $\sigma >1$.
We find that there are no zeroes for each of these values of $N$. Combining this with Table \ref{table1}, one proves the following  theorem.

\begin{theorem}\label{MainT}
For $1\leq N \leq 18$ and $N=20, 21, 28$ there are no zeroes of $\zeta_{N}(s)$ in the region $\sigma>1$; for all other positive $N$ there are infinitely many such zeroes.
\end{theorem}

\section{Numerical Computation}\label{comp}

\subsection{Interval Arithmetic}

Almost all real numbers are not exactly representable by any finite-precision, floating-point system such as the $64$ bit IEEE implementation available on most modern processors. Thus any computation involving such a floating point system will, unless we are very lucky, only produce an approximation to the true result. One way of managing this is to use interval arithmetic (see, for example, \cite{Moore1966} for a good introduction). Instead of storing a floating point number that is an approximation to the value we want, we store an interval bracketed by two floating point numbers that contains the true value. 
%We then define operators $+$, $-$, $\times$ and $\div$ together with functions such as $\sqrt{}$ and $\exp$ to handle such intervals faithfully.

Interval arithmetic has been used to manage the accumulation of round-off and truncation errors. In this paper, we exploit the technique to get zero free regions  rigorously. As an example, consider the function $f:\mathbb{R}\rightarrow\mathbb{R}$ defined by
\begin{equation*}
f(x)=x^2-4x+3.
\end{equation*}
Suppose we wish to demonstrate that $f$ has no zeroes for $x\in[4,5]$. Then we can compute
\begin{equation*}
f([4,5])=[16,25]-[16,20]+3=[-1,12].
\end{equation*}
Since this is inconclusive, we try again, but this time with the interval split in two. We have
\begin{equation*}
f([4,4.5])=[16,20.25]-[16,18]+3=[1,7.25]
\end{equation*}
and
\begin{equation*}
f([4.5,5])=[20.25,25]-[18,20]+3=[3.25,10]
\end{equation*}
and we have our result\footnote{Note that if we had written $f(x)=(x-1)(x-3)$ then $f[4,5]=[3,4]\cdot[1,2]=[3,8]$ which is the ``correct'' result. This sensitivity is common in expressions involving intervals.}.

\subsection{Description of algorithm}
We first note that we need not search in all of $\sigma>1$ to find zeroes of $\zeta_{N}(s)$. Spira \cite[Thm.\ 1]{SpiraSections1} proved that all zeroes of $\zeta_{N}(s)$ must have real part less than $1.85$; this was sharpened in \cite[Theorem 3.1]{BorweinZero} to $1.73$. We therefore need only consider $\sigma\in(1, 1.73]$. We can improve this for some values of $N$, but, as we shall see in \S \ref{2.3}, this is more than sufficient for our purposes.

Let us consider the case $N=28$. Let $p$ denote a prime and let $\theta_p=t\log p$. Hence we have
\begin{equation*}
\zeta(\sigma+it)=1+\frac{\exp(-i\theta_2)}{2^\sigma}+\frac{\exp(-i\theta_3)}{3^\sigma}+\frac{\exp(-i 2\theta_2)}{4^\sigma}+\ldots+\frac{\exp(-i (2\theta_2+\theta_7))}{28^\sigma}
\end{equation*}
and we will now write $\zeta_{28}(\sigma,\theta_2,\ldots,\theta_{23})$ for $\zeta_{28}(\sigma+it)$ under such a change of variables.

It would appear that we need to examine the space $\sigma\in(1, 1.73]$,
 $\theta_p\in[0,2\pi)$ for $p\leq 23$, for zeroes. In fact we can do considerably better. First, we observe that $\theta_{17}$, $\theta_{19}$ and $\theta_{23}$ only appear once in the sum. Call the sum without those three terms $\zeta_{28'}(\sigma,\theta_2,\ldots,\theta_{13})$. Then $\zeta_{28}$ cannot have a zero if there is no $\sigma,\theta_2,\ldots,\theta_{13}$ such that
\begin{equation*}
\left|\zeta_{28'}(\sigma,\theta_2,\ldots,\theta_{13})\right|\leq 17^{-\sigma}+19^{-\sigma}+23^{-\sigma}.
\end{equation*}

We can go further. The $\theta_{11}$ and $\theta_{13}$ terms only appear on their own or in conjunction with $\theta_2$. We write $a=11^{-\sigma}$, $b=22^{-\sigma}$, $c=13^{-\sigma}$ and $d=26^{-\sigma}$. Then a little high school geometry (the cosine rule to be precise) tells us that 
\begin{equation*}
\left|\frac{\exp(-i\theta_{11})}{11^\sigma}+\frac{\exp(-i(\theta_{11}+\theta_{2}))}{22^{-\sigma}}\right|\leq\sqrt{a^2+b^2+2ab\cos \theta_2},
\end{equation*}
and
\begin{equation*}
\left|\frac{\exp(-i\theta_{13})}{13^\sigma}+\frac{\exp(-i(\theta_{13}+\theta_{2}))}{26^{-\sigma}}\right|\leq\sqrt{c^2+d^2+2cd\cos \theta_2}.
\end{equation*}

Call $\zeta_{28''}(\sigma,\theta_2,\theta_3,\theta_5,\theta_7)$ the result obtained by removing the $n=11,13,22$ and $26$ terms from $\zeta_{28'}$. With $a,b,c$ and $d$ as above, define
\begin{equation*}
f(\sigma,\theta_2)=17^{-\sigma}+19^{-\sigma}+23^{-\sigma}+\sqrt{a^2+b^2+2ab\cos\theta_2}+\sqrt{c^2+d^2+2cd\cos\theta_2}.
\end{equation*}
Then $\zeta_{28}$ cannot have a zero if there is no $\sigma,\theta_2,\ldots,\theta_{7}$ such that
\begin{equation*}
\left|\zeta_{28''}(\sigma,\theta_2,\theta_3,\theta_5,\theta_7)\right|\leq f(\sigma,\theta_2).
\end{equation*}

Our algorithm is as follows. Divide $\sigma$, $\theta_2$, $\theta_3$, $\theta_5$ and $\theta_7$ into small intervals that cover $[1,1.73]$ and $[0,2\pi]^4$ respectively. We refer to any choice of five such intervals as a ``box''. Push all possible boxes onto the stack. While the stack is not empty, pop off a box and compute an interval $z$ containing $|\zeta_{28''}|$ for that box. Compute an interval containing $f(\sigma,\theta_2)$. If the interval $z-f(\sigma,\theta_2)$ is wholly positive, then that box did not contain any zeroes, so discard it. If the interval is wholly negative, then terminate with failure.\footnote{We believe that this condition indicates the presence of infinitely many zeroes. We are grateful to the referee for suggesting a means by which one might seek to establish this, based on \cite{AvellarandHale},\cite{Dubonetal} and \cite{SepulcreandVidal}. However, the weaker statement is sufficient for our purposes and we do not pursue this line of thought further.} If the interval straddles zero, then divide the box into $16$ smaller boxes by halving the intervals for the $\theta_p$, and push these new boxes onto the stack.

\subsection{Details of the implementation}\label{2.3}

We implemented this algorithm in ``C++'' using our own double precision interval package written in assembler. This exploits an idea of Lambov \cite{Lambov2008} to make efficient use of the SSE instruction set of modern processors and uses CRMLIB \cite{Muller2010} to implement the transcendental functions.

We divided the interval for $\sigma$ into $16$ sub-intervals $[1+(2^n-1)\cdot 2^{-16},1+(2^{n+1}-1)\cdot2^{-16}]$ for $0\leq n \leq 15$. Therefore, the first interval checked was $\sigma=[1,1+2^{-16}]$ and the last\footnote{Note that this covered a wider interval than was strictly necessary.} was $\sigma=[\frac{3}{2}-2^{-16},2-2^{-16}]$. Each of these intervals for $\sigma$ was handled by a single core of a compute node of the University of Bristol's Bluecrystal Phase III cluster \cite{ACRC2015}\footnote{A single node of Phase III contains two 8-core Intel E5-2670 Sandybridge processors running at $2.6$GHz.}. Within a single core, the intervals for $\theta_2,\theta_3,\theta_5$ and $\theta_7$ were initially divided into $16$, $8$, $4$ and $2$ sub-intervals respectively, for a total of $1,024$ boxes. Since $\theta_2$ contributed to more terms than the other variables, it made sense to start with a narrower search here: this seemed to work well in practice.

Table \ref{table2} shows the data for $N=28$ and $\sigma\in[1,1+2^{-16}]$. At each iteration, a box could result in $16$ new boxes; at first this is what we see. We see that after the second iteration, the remaining search space decreases dramatically.

\begin{table}[ht]
\caption{Number of boxes at each iteration for $N=28$ and $\sigma\in[1,1+2^{-16}]$}
\label{table2}
\centering
\begin{tabular}{c c c}
\hline\hline
Iteration & Number of Boxes &  Coverage (\%)\\[0.5ex]\hline
$1$ & $1,024$ & $100$\\
$2$ & $16,256$ & $99.2$\\
$3$ & $45,920$ & $17.5$\\
$4$ & $118,560$ & $2.83$\\
$5$ & $170,048$ & $0.25$\\
$6$ & $195,920$ & $0.018$\\
$7$ & $212,960$ & $0.0012$\\
$8$ & $82,016$ & $0.000030$\\
    \hline\hline
  \end{tabular}
\end{table}

We ran this algorithm for 
%$N=\{10,11,12,13,14,15,16,17,18,20,21,28\}$ 
those $N$ in (\ref{big}) and in every case confirmed that $\zeta_N(s)$ has no zeroes for $\sigma\geq 1$. Checking each $N$ took much less than a minute elapsed time using $16$ cores, with $N=21$ taking the longest at $30$ seconds.

\end{document}